\documentclass[11pt]{article}
\usepackage{amssymb,amsfonts,amsmath,latexsym,epsf,tikz,url}

\newtheorem{theorem}{Theorem}[section]
\newtheorem{proposition}[theorem]{Proposition}

\newtheorem{corollary}[theorem]{Corollary}

\newtheorem{remark}[theorem]{Remark}
\newtheorem{definition}[theorem]{Definition}

\newcommand{\proof}{\noindent{\bf Proof.\ }}
\newcommand{\qed}{\hfill $\square$\medskip}

\begin{document}

\title{ Introduction to total dominator edge chromatic  number }

\author{
Nima Ghanbari \and Saeid Alikhani$^{}$\footnote{Corresponding author}
}

\date{\today}

\maketitle

\begin{center}
Department of Mathematics, Yazd University, 89195-741, Yazd, Iran\\
{\tt n.ghanbari.math@gmail.com, alikhani@yazd.ac.ir}
\end{center}


\begin{abstract}
	We introduce the total dominator edge chromatic number of a graph $G$. A total dominator edge coloring (briefly TDE-coloring)  of  $G$ is a proper edge coloring of $G$ in which each edge of the graph is adjacent to every edge of some color class. The total dominator edge chromatic number (briefly TDEC-number)  $\chi'^t_d(G)$ of $G$ is the minimum number of color classes in a TDE-coloring of $G$. We obtain some properties of $\chi'^t_d(G)$ and  compute this parameter for specific graphs. We   examine the effects on $\chi'^t_d(G)$ when $G$ is modified by operations on vertex and  edge of $G$. Finally, we consider the $k$-subdivison of $G$ and study TDEC-number of this kind of graphs.  
\end{abstract}

\medskip

\noindent{\bf AMS Subj.\ Class.}: 05C15, 05C69


\section{Introduction}


Let $G=(V,E)$ be a simple  graph and $k \in \mathbb{N}$. A mapping $f : V (G)\longrightarrow \{1, 2,...,k\}$ is
called a $k$-proper  coloring of $G$, if $f(u) \neq f(v)$ whenever the vertices $u$ and $v$ are adjacent
in $G$. A color class of this coloring, is a set consisting of all those vertices
assigned the same color. If $f$ is a proper coloring of $G$ with the coloring classes $V_1, V_2,..., V_{k}$ such
that every vertex in $V_i$ has color $i$, then sometimes write simply $f = (V_1,V_2,...,V_{k})$.  The chromatic number $\chi(G)$ of $G$ is
the minimum number of colors needed in a proper coloring of a graph.

The total dominator coloring, abbreviated TD-coloring studied in \cite{Adel,Adel2,Vij1,Vij2}. Let $G$ be a graph with no
isolated vertex, the total dominator coloring  is a proper coloring of $G$ in which each vertex of the graph is adjacent
to every vertex of some (other) color class. The total dominator chromatic number, abbreviated TDC-number, $\chi_d^t(G)$ of $G$ is the minimum number of color classes in a TD-coloring of $G$.   Computation of the TDC-number is NP-complete (\cite{Adel}).
The TDC-number of some graphs has computed in
\cite{Adel}. Also Henning in \cite{GCOM} established the  lower and the upper bounds on the TDC-number
of a graph in terms of its total domination number $\gamma_t(G)$. He has shown that,  every
graph $G$ with no isolated vertex satisfies $\gamma_t(G) \leq \chi_d^t (G)\leq \gamma_t(G) + \chi(G)$.
The properties of TD-colorings in trees has studied in \cite{GCOM,Adel}. Trees $T$ with $\gamma_t(T) =\chi_d^t(T)$ has characterized
in \cite{GCOM}. We have examined the effects on $\chi_d^t(G)$ when $G$ is modified by operations on the vertex and the edge of $G$, and  the TDC-number of some operations on two graphs studied in \cite{nima2}. 

\medskip
Motivated by TDC-number of a graph, we consider the proper edge coloring of $G$ and introduce the total dominator edge chromatic number (TDEC-number) of $G$, $\chi'^t_d(G)$, obtain some properties of $\chi'^t_d(G)$ and   compute this parameter for specific graphs, in the next section. In Section 3, we  examine the effects on $\chi'^t_d(G)$ when $G$ is modified by operations on vertex and  edge of $G$. Finally in Section 4, we study the TDEC-number of $k$-subdivision of graphs.

\section{Introduction to total dominator edge chromatic number}

In this section, we state the definition of total dominator edge chromatic number and obtain this parameter for some
specific graphs. 

\begin{definition}
	A {\rm total dominator edge coloring}, briefly TDE-coloring, of a graph $G$ is a proper edge coloring of $G$ in which each edge of the graph is adjacent to every edge of some color class. The total dominator edge chromatic number (TDEC-number) $\chi'^t_d(G)$ of $G$ is the minimum number of color classes in a TDE-coloring  of $G$. A $\chi'^t_d(G)$-coloring of $G$ is any total dominator edge coloring of $G$ with $\chi'^t_d(G)$ colors.
\end{definition}

\begin{remark} \label{1} 
For every graph $G$ with maximum degree $\Delta(G)$, $\chi'^t_d (G) \geq \Delta(G).$ This inequality is sharp. As an example, for the star graph $K_{1,n}$, $\chi'^t_d (K_{1,n}) =n.$
\end{remark}

The following theorem gives the total dominator edge chromatic number of a path.

\begin{theorem}\label{Pn}
If $P_n$ is  the path graph of order $n\geq 9$, then
\[
\chi'^t_d(P_n)=\left\{
	\begin{array}{ll}
	{\displaystyle
		2k+2}&
	\quad\mbox{if $n=4k+1$,}\\[15pt]
{\displaystyle
		2k+3}&
	\quad\mbox{if $n=4k+2$,}\\[15pt]
{\displaystyle
		2k+4}&
	\quad\mbox{if $n=4k+3$, $4k+4$.}
	\end{array}
	\right.
	\]
Also $ \chi'^t_d(P_3)= \chi'^t_d(P_4)=2$, $ \chi'^t_d(P_5)=3$, $ \chi'^t_d(P_6)=\chi'^t_d(P_7)=4$ and $ \chi'^t_d(P_8)=5$.
	
\end{theorem}

\begin{figure}
	\begin{center}
		\includegraphics[width=2.3in]{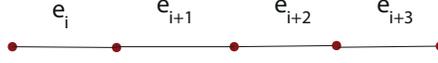}
		\caption{Four consecutive edges. }
		\label{Four}
	\end{center}
\end{figure}

\proof 
It is easy to show that $\chi'^t_d (P_3)=\chi'^t_d (P_4)=2$, $\chi'^t_d (P_5)=3$, $\chi'^t_d (P_6)=\chi'^t_d (P_7)=4$ and $\chi'^t_d (P_8)=5$. Suppose that  $n\geq 9$. 
First we show that in a TDE-coloring,  for each four consecutive edges we shall  use at least two  new colors. We consider two cases.  If an used  color assign  to edge  $e_{i+1}$, then we need to assign a new color to  the edge  $e_{i+2}$ and $e_{i+3}$ to have a TDE-coloring (see Figure \ref{Four}). If a new color assign to the edge $e_{i+1}$, then we have to assign  a new color to $e_{i+2}$ or $e_{i}$ to have a TDE-coloring. So we need at least two new colors in every four consecutive vertices. 

If  $n=4k+1$, for some $k\in \mathbb{N}$, then we give a TDE-coloring  for $P_{4k+1}$ which use only two new colors in every four consecutive edges. Define a function $f_0$ on $E(P_{4k})$ such that for any edge $e_i$,

  \[
   f_0(e_i)=\left\{
	\begin{array}{ll}
	{\displaystyle
		1} &
	\quad\mbox{if $i=1+4s$,}\\[15pt]
	{\displaystyle
		2} &
	\quad\mbox{if $i=4s$.}\\[15pt]
	{\displaystyle
		2} &
	\quad\mbox{if $i=4s$.}
	\end{array}
	\right.
	\]
where $s$ is a natural number and for any $e_i$, $i \neq 4s$ and $i\neq 4s+1$, $f_0(e_i)$ is a new number. Then this coloring is a TDE-coloring  of $P_{4k+1}$ 
with the minimum number $2k+2$ colors.

If  $n=4k+2$, for some $k\in \mathbb{N}$, then  we first color the $4k-4$ edges using  $f_0$. Now for the rest of edges we define $f_1$ as  $f_1(e_{4k-3})=1$, $f_1(e_{4k-2})=2k+1$, $f_1(e_{4k-1})=2k+2$, $f_1(e_{4k})=2k+3$ and $f_1(e_{4k+1})=2$. Since  for every five consecutive edges we have to use at least three  new colors, so this  edge coloring is  a TDE-coloring  of $P_{4k+1}$ with the minimum number $2k+3$ colors.

If  $n=4k+3$, for some $k\in \mathbb{N}$, then using $f_0$  we color the $4k-4$ edges. Now for the rest of edges define $f_2$ as $f_2(e_{4k-3})=1$, $f_2(e_{4k-2})=2k+1$, $f_2(e_{4k-1})=2k+2$, $f_2(e_{4k})=2k+3$, $f_2(e_{4k+1})=2k+4$ and $f_2(e_{4k+2})=2$. Since  for every six consecutive edges we have to use at least four  new colors, so this edge coloring  is a TDE-coloring  of $P_{4k+2}$ with the minimum number $2k+4$ colors.

If  $n=4k+4$, for some $k\in \mathbb{N}$, then using $f_0$  we color  the $4k-4$ edges and  for the rest of edges define $f_3$ as $f_3(e_{4k-3})=1$, $f_3(e_{4k-2})=2k+1$, $f_3(e_{4k-1})=2k+2$, $f_3(e_{4k})=2$, $f_3(e_{4k+1})=2k+3$, $f_3(e_{4k+2})=2k+4$ and $f_3(e_{4k+2})=2$. This coloring is a TDE-coloring  of $P_{4k+2}$ with the minimum number $2k+4$ colors. So we have the result. \qed
												

\begin{figure}
\begin{center}
\includegraphics[width=4in]{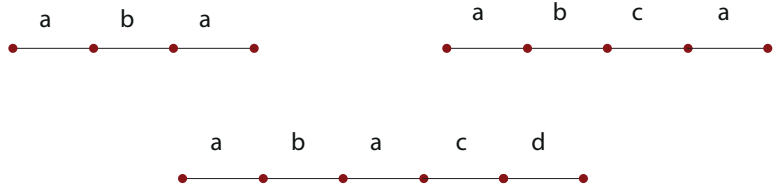}
\caption{TDE-coloring of $P_4$, $P_5$ and $P_6$.}
\label{P456}
\end{center}
\end{figure}

						
\begin{figure}
\begin{center}
\includegraphics[width=1.5in]{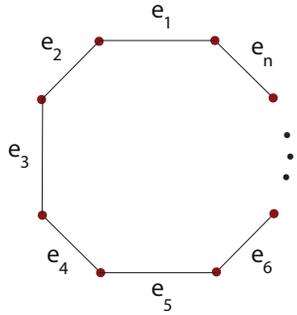}
\caption{Cycle graph of order $n$, $C_n$.}
\label{Cne}
\end{center}
\end{figure}

\begin{theorem}\label{Cn}
If $C_n$ is the cycle graph of order $n\geq 8$, then
\[
      \chi'^t_d(C_n)=\left\{
	\begin{array}{ll}
	{\displaystyle
		2k+2,}&
	\quad\mbox{if $n=4k$,}\\[15pt]
{\displaystyle
		2k+3,}&
	\quad\mbox{if $n=4k+1$,}\\[15pt]
{\displaystyle
		2k+4,}&
	\quad\mbox{if $n=4k+2$, $4k+3$. }
	\end{array}
	\right.
	\]	
Also $\chi'^t_d (C_3)=3$, $\chi'^t_d (C_4)=2$, $\chi'^t_d (C_5)=\chi'^t_d (C_6)=4$ and $\chi'^t_d (C_7)=5$.
\end{theorem}

\proof 
It is similar to the Proof of Theorem \ref{Pn}. \qed

The following corollry is an immediate consequence of Theorems \ref{Pn} and \ref{Cn}.

\begin{corollary}
For every $n\geq6$, $\chi'^t_d(P_n)=\chi'^t_d(C_{n-1})$.
\end{corollary}

			The following theorem present  a lower bound for TDEC-number of graphs $G$ which have the graph  $P_6$ as induced subgraph. 	
				
\begin{theorem}\label{delta1}
					If  $G$ is a connected graph containing  $P_6$ as an induced subgraph, then $\chi'^t_d(G) \geq \Delta (G)+2$.  More generally, if the path graph $P_n$ is an induced subgraph of $G$, then $\chi'^t_d(G) \geq \Delta (G)+\chi'^t_d(P_{n-2})$.
				\end{theorem}
				\proof
				We assign  $\Delta (G)$ colors to the edges which are incident  to the vertex with maximum degree $\Delta(G)$. Now we consider $P_6$ as  induced subgraph of $G$. 
				As we have seen in the Proof of Theorem \ref{Pn}, we need at least two new colors for each four consecutive edges. So we have $\chi'^t_d(G) \geq \Delta (G)+2$. The proof of inequality  $\chi'^t_d(G) \geq \Delta (G)+\chi'^t_d(P_{n-2})$ is similar. \qed

\begin{figure}
					\begin{center}
						\includegraphics[width=2.5in]{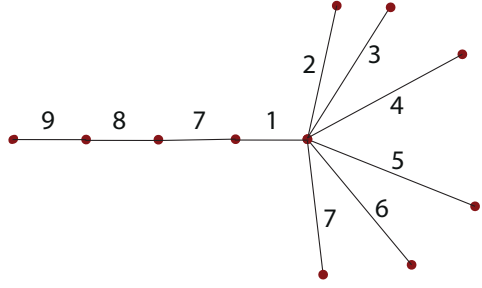}
						\caption{Graph $G$ with $\chi'^t_d(G)= \Delta (G)+2$.}
						\label{G}
					\end{center}
				\end{figure}

				\medskip 
				\begin{remark}
					The graph $G$ in Figure \ref{G} and its coloring shows that the lower bound in Theorem \ref{delta1} is sharp.
				\end{remark} 												
				
				\begin{theorem}
					For every $n\in \mathbb{N}$, 	 $2n-1\leq \chi'^t_d(K_{2n})\leq4n-2$ and $2n\leq \chi'^t_d(K_{2n+1})\leq 4n-1$.
				\end{theorem}

				\proof
				 The lower bounds follow from Remark \ref{1}. To obtain the upper bound, suppose that $V(K_{2n+1})=\lbrace u_1,\ldots,u_{2n+1} \rbrace$. By removing the vertex $u_{2n+1}$, we have the complete graph $K_{2n}$. We know that  $\chi'(K_{2n})=2n-1$. So we color the edges of  $K_{2n}$ with $2n-1$ colors. Now we add the vertex $u_{2n+1}$ and make $K_{2n+1}$ and assign  the new colors $2n,2n+1,\ldots,4n-1$ to new  edges. This is a TDE-coloring  for $K_{2n+1}$. Therefore $ \chi'^t_d(K_{2n+1})\leq 4n-1$. By the similar method  we have $ \chi'^t_d(K_{2n})\leq 4n-2$.
				\qed


				\begin{theorem}\label{bipartite}
	\begin{enumerate}
		\item [(i)]
 For every  $n\neq m$,  $max \lbrace n,m\rbrace \leq \chi'^t_d(K_{n,m})\leq m+n-1$. 
 \item[(ii)] 
  For every $n\in \mathbb{N}$, $ n \leq \chi'^t_d(K_{n,n})\leq 2n$. 
  \end{enumerate} 
				\end{theorem}

					\proof
						\begin{enumerate}
							\item [(i)]
			The lower bounds follow from Remark \ref{1}. To obtain the upper bound, suppose that $V(K_{n,m})=X\cup Y$, where $X=\{u_1,\ldots,u_m\}$ and 
			$Y=\{u_{m+1},\ldots , u_{m+n}\}$ and $m\geq n$. We have the following cases:
						
			\noindent Case 1) $m=n+1$. By removing the vertex $u_1$, we have the complete bipartite graph $K_{n,n}$. We know that  $\chi'(K_{n,n})=n$. So we color the edges of  $K_{n,n}$ with $n$ colors. Now we add the vertex $u_1$ and make $K_{n+1,n}$ and assign the new colors $n+1,n+2,\ldots,2n$ to new  edges. This is a TDE-coloring  for $K_{n,m}$ and we have $\chi'^t_d(K_{n,m})\leq 2n = m+n-1$.
			
		\noindent 	Case 2) $m>n+1$. By removing the vertex $u_1$, we have the complete bipartite graph $K_{m-1,n}$. We know that  $\chi'(K_{m-1,n})=m-1$. So we color the edges of  $K_{m-1,n}$ with $m-1$ colors. Now we add the vertex $u_1$ and make $K_{m,n}$ and assign  the new colors $m,m+1,\ldots,m+n-1$ to new  edges. This is a TDE-coloring  for $K_{m,n}$ and we have $\chi'^t_d(K_{m,n})\leq m+n-1$.
			
	\item[(ii)] 	In this part we have $m=n$.  By removing the vertex $u_1$, we have the complete bipartite graph $K_{n-1,n}$. We know that  $\chi'(K_{n-1,n})=n$. So we color the edges of  $K_{n-1,n}$ with $n$ colors. Now we add the vertex $u_1$ and make $K_{n,n}$ and assign the new colors $n+1,n+2,\ldots,2n$ to new  edges. This is a TDE-coloring  for $K_{n,n}$ and we have $\chi'^t_d(K_{n,n})\leq 2n $.	
				\qed			
	 \end{enumerate}

\medskip 
\begin{remark}
	 The lower bounds in parts (i) and (ii) of Theorem \ref{bipartite} are sharp. It suffices to consider $K_{3,2}$ and  $K_{2,2}=C_4$, respectively.  Note that   $\chi'^t_d(K_{3,2})=3$ and $\chi'^t_d(C_4)=2$. Also the upper bound of part (i) is sharp. It suffices to consider the star graph $K_{1,6}$. Note that  $\chi'^t_d(K_{1,6})=1+6-1=6$.
	 \end{remark}

\begin{figure}[ht]
		\hspace{2cm}
		\begin{minipage}{4.5cm}
			\includegraphics[width=\textwidth]{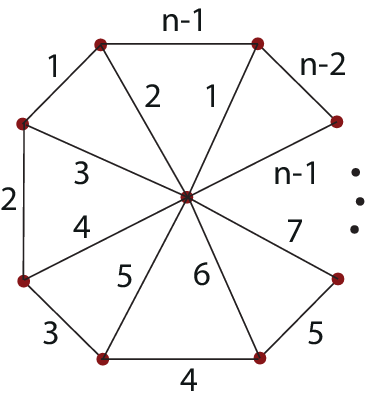}
					\end{minipage}
					\hspace{1cm}
		\begin{minipage}{4.5cm}
			\includegraphics[width=\textwidth]{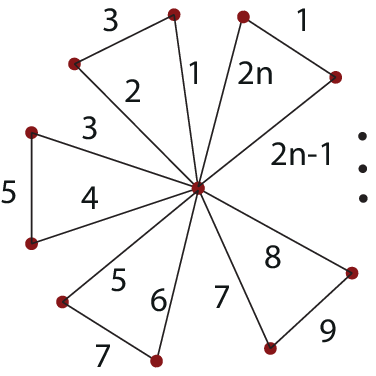}
		\end{minipage}
		\caption{TDE-coloring of wheel and friendship graph of order $n$}
		\label{Wne}
	\end{figure}

		 Let $n$ be any positive integer and $F_n$ be the friendship graph with $2n+1$ vertices and $3n$ edges, formed by the join of $K_1$  with $nK_2$. 
	By Remark \ref{1},  and  TDE-coloring which  has shown in  Figure \ref{Wne},   we have the following result for the wheel of order $n$, $W_n$ and the friendship graph $F_n$.											
												
\begin{theorem}\label{Wn}
\begin{enumerate}
\item[(i)] 
For any $n\geq 3$, 									
$\chi'^t_d (W_n)=n-1.$

\item[(ii)] 
For  $n\geq 2$,  $\chi'^t_d (F_n)=2n.$
\end{enumerate}					
\end{theorem}

\section{TDEC-number of some operations on a graph}
												
The graph $G-v$ is a graph that is made  by deleting the vertex $v$ and all edges incident  to $v$ from the graph $G$ and the graph $G-e$ is a graph that obtained from $G$ by simply removing the  edge $e$. In this section we  present   bounds for TDEC-number  of $G-v$ and $G-e$.	We begin with $G-e$. 											
												
\begin{theorem}\label{removal} 
If $G$ is  a connected graph, and $e\in E(G)$ is not a bridge of $G$, then 
$$ \chi'^t_d (G)-2\leq \chi'^t_d (G-e)\leq  \chi'^t_d (G)+2 .$$
\end{theorem}
												
\begin{figure}
\begin{center}
\includegraphics[width=6in]{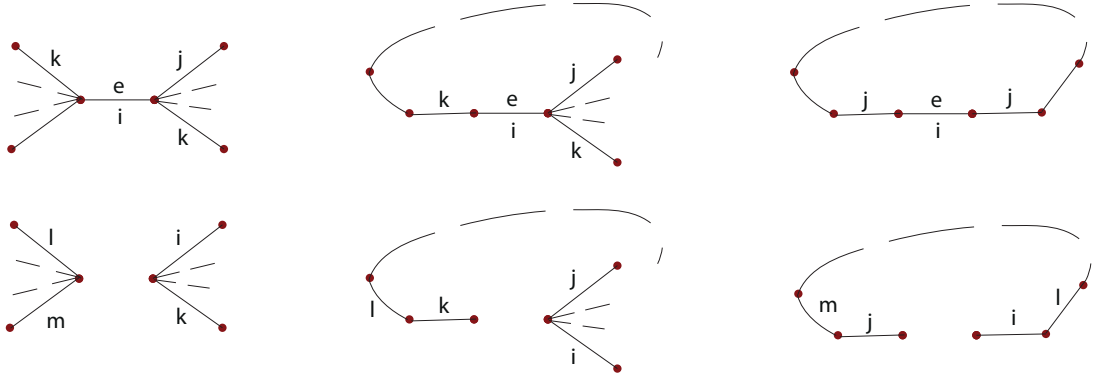}
\caption{Cases which has considered in the proof of Theorem \ref{removal}.}
	\label{G-e1}
\end{center}	
\end{figure}

\begin{figure}
\begin{center}
\includegraphics[width=4.5in]{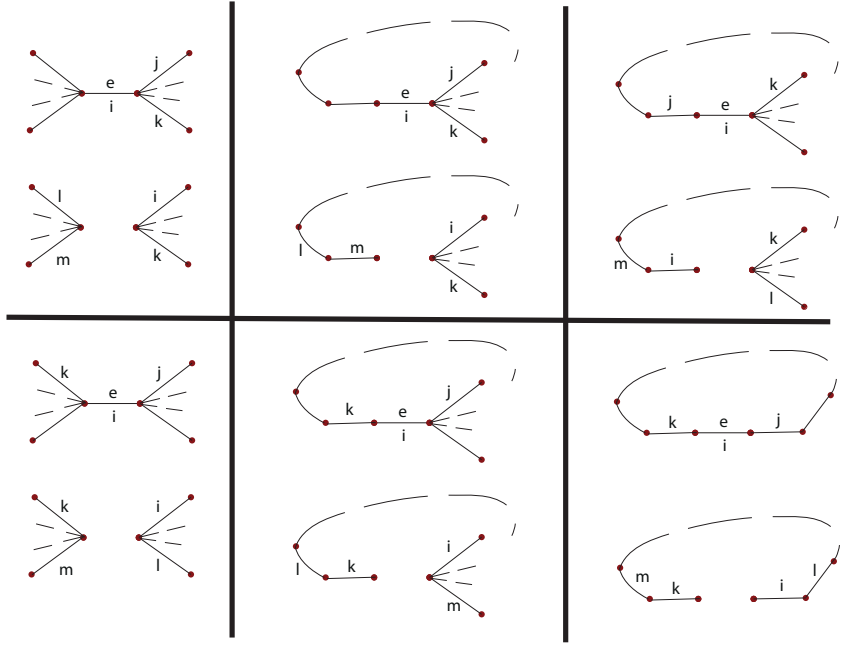}
\caption{Another cases which has considered in the proof of Theorem \ref{removal}.}
\label{G-e2}
\end{center}\end{figure}

\proof 
First we prove the right inequality.  Suppose that the edge $e$ in a TDE-coloring of $G$ has color $i$. If no edges of $G$ use the color class $i$, then TDE-coloring of $G$ is  a TDE-coloring of $G-e$, too. So $\chi'^t_d (G-e)\leq  \chi'^t_d (G)$. If some edges of $G$ use the color class $i$ in TDE-coloring, then we have at most two edges with color $i$.
If  two edges of $G$ have color $i$, then removing $e$ does not effect on TDE-coloring  and any edge uses the old color class in TDE-coloring of $G$. So $\chi'^t_d (G-e)\leq  \chi'^t_d (G)$. If  only one edge $e$ has the  color $i$, then we change the color of some edges in $G-e$ to have a TDE-coloring for $G-e$. In this case the edge $e$ uses some color class, say $k$, and is adjacent to all color class $k$. We can not have more than two $k$ in this case. Suppose that we have two $k$. Then we have only three cases for the graph $G$ as we see in Figure \ref{G-e1}. In Figure \ref{G-e1}, the colors $l$ and $m$ are new colors and we only change the color of some edges in $G-e$ and assign  the other edges their old color in $G$. This coloring is  a TDE-coloring  for $G-e$. In any case we do not use more than two new colors. Therefore we have $\chi'^t_d (G-e)\leq  \chi'^t_d (G)+2$. 
Now suppose that we have only one color $k$. Then we have only six cases for the graph $G$ as we see in Figure \ref{G-e2}. In Figure \ref{G-e2}, the colors $l$ and $m$ are new colors and we only change the color of some edges in $G-e$ and assig the other edges their old color in $G$. This  kind of coloring is  a TDE-coloring  for $G-e$. In any case we do not  use more than two new colors. 												
So we have $\chi'^t_d (G-e)\leq  \chi'^t_d (G)+2$.
												
Now we prove that $ \chi'^t_d (G)-2\leq \chi'^t_d (G-e)$. To do this,  first we color $G-e$ and then we add edge $e$. 
We assign new color $i$ to $e$ and new color $j$ to one edge which is adjacent to $e$. So we have a TDE-coloring for $G$ and $ \chi'^t_d (G)\leq \chi'^t_d (G-e)+2$. Therefore we have the result.\qed

\begin{theorem}
If $G$ is a connected graph, and  $v\in V(G)$ is not a cut vertex of $G$, then 
$$\chi'^t_d (G)-deg(v)\leq \chi'^t_d (G-v) \leq \chi'^t_d (G)+deg(v).$$
\end{theorem}

\proof
First we prove the left  inequality. We give a TDE-coloring to $G-v$,  add $v$ and all the corresponding edges. Then we assign  $deg(v)$ new colors to these 
edges and do not change the color of other edges. So this is a TDE-coloring of $G$ and $\chi'^t_d (G)\leq \chi'^t_d (G-v) + deg(v)$. 
												
For the right   inequality, first we give a TDE-coloring  to $G$. In this case, since $v$ is not a cut vertex, each edge which is adjacent to an edge with endpoint $v$ has an other adjacent edge too. we change the color of this edge to a new color and do this $deg(v)$ times and do not change the color of the other edges.  So this is a TDE-coloring of $G-v$ and $\chi'^t_d (G-v) \leq \chi'^t_d (G)+deg(v)$. Therefore we have the result. 	\qed

The following theorem is an immediate consequence of  Theorems \ref{Cn} and \ref{Wn}.											
\begin{theorem}
There is a connected graph $G$, and a vertex  $v\in V(G)$ which is not a cut vertex of $G$ such that 	 $\vert \chi'^t_d(G)  - \chi'^t_d (G-v) \vert $ can be arbitrarily large.
\end{theorem}
								
											
In a graph $G$, contraction of an edge $e$ with endpoints $u,v$ is the replacement of $u$ and $v$ with a single vertex such that edges incident to the new vertex are the edges other than $e$ that were incident with $u$ or $v$. The resulting graph $G/e$ has one less edge than $G$. We denote this graph by $G/e$. We end this section with the following theorem which
gives bounds for $\chi'^t_d(G/e)$.

\begin{theorem}\label{goe}
If  $G$ is a connected graph and $e=uv\in E(G)$, then
$$\chi'^t_d(G)-2 \leq \chi'^t_d(G/e) \leq \chi'^t_d(G)+min\{deg(u),deg(v)\}-1 .$$
\end{theorem}
\proof
First we prove the left inequality. We give a TDE-coloring to $G/e$,  add $e$  and assign it a new color, say  $i$ and change the color of one of its adjacent edges to new color $j$ and do not change other colors. This is a TDE-coloring  of $G$. So we have $\chi'^t_d(G) \leq \chi'^t_d(G/e)+2$. 
For the right inequality, we give a TDE-coloring  to $G$. Suppose that $min\{deg(u),deg(v)\}=deg(u)$. Now we make $G/e$ and change the color of adjacent edges of  $e$ with the endpoint $u$ to new colors. So we have the result.			\qed

	\medskip 
\begin{remark}
	 The lower bound in Theorem \ref{goe} is sharp. It suffices to consider the cycle graph $C_5$ as $G$. Note that  $\chi'^t_d(C_5)=4$ and   $\chi'^t_d(C_4)=2$. 
	 \end{remark}

 \section{TDEC-number of $k$-subdivision of a graph}

 The \textit{ $k$-subdivision} of $G$, denoted by $G^{\frac{1}{k}}$, is constructed by replacing each edge $v_iv_j$ of $G$ with a path of length $k$., say $P_{\{v_i,v_j\}}$. These $k$-paths are called \textit{superedges}, any new vertex is an internal vertex, and is denoted by $x^{\{v_i,v_j\}}_l$ if it belongs to the superedge $P_{\{v_i,v_j\}}$, $i<j$ with  distance $l$ from the vertex $v_i$, where $l \in \{1, 2, \ldots , k-1\}$.  Note that for $k = 1$, we have $G^{1/1}= G^1 = G$, and if the graph $G$ has $v$ vertices and $e$ edges, then the graph $G^{\frac{1}{k}}$ has $v+(k-1)e$ vertices and $ke$ edges. The total dominator chromatic number of a graph has studied in \cite{subdivision}.  In this section we study TDEC-number of $k$-subdivision of a graph. In particular, we obtain some  bounds for $\chi'^t_d(G^{\frac{1}{k}})$ and prove that for any $k\geq 2$, $\chi'^t_d(G^{\frac{1}{k}}) \leq \chi'^t_d(G^{\frac{1}{k+1}})$.

  \begin{theorem}
 	If $G$ is  a graph with $m$ edges, then $\chi'^t_d(G^{\frac{1}{k}}) \geq m$, for $k\geq 3$. 
 \end{theorem}
  \proof
 For $k=3$, in any superedge $P^{\{v,w\}}$ such as $\{v, x^{\{v,w\}}_1,  x^{\{v,w\}}_2,   w\}$ ,The edge $x^{\{v,w\}}_1 x^{\{v,w\}}_{2}$ need to use a new color in at least one of its adjacent edges, and we cannot use this color in any other superedges. 
 So we have the result. \qed

  \begin{theorem}\label{ubboun2p}
  	If  $G$ is a connected graph with $m$ edges and $k\geq 2$, then 
  	\begin{equation*}
  	\chi'^t_d(P_{k+1}) \leq \chi'^t_d(G^{\frac{1}{k}})\leq m \chi'^t_d(P_{k+1}).
  	\end{equation*}
  \end{theorem}
  \proof  First we prove the the right inequality. Suppose that  $e=uu_1$ be an arbitrary edge of $G$. This edge is replaced with   the super edge $P^{\{u,u_1\}}$  in $G^{\frac{1}{k}}$, with vertices $\{u, x^{\{u,u_1\}}_1, \ldots , x^{\{u,u_1\}}_{k-1}, u_1\}$. We color this superedge with  $\chi'^t_d(P_{k+1})$ colors as a total dominator edge coloring of $P_{k+1}$. We do this for all superedges. 
  Thus we need  at most $m \chi'^t_d(P_{k+1})$ new colors for a total dominator edge coloring of $G^{\frac{1}{k}}$.
  
  For the left inequality, if $G$ is a path then the result is true. So we suppose that $G$ is a connected graph which is not a path. Let  $P^{\{v,w\}}$ be an arbitrary superedge of $G^{\frac{1}{k}}$ with vertex set  $\{v, x^{\{v,w\}}_1, \ldots , x^{\{v,w\}}_{k-1}, w\}$. Since $G$ is not a path, so at least one of $v$ and $w$ is adjacent to some vertices of $G^{\frac{1}{k}}$ except $x^{\{v,w\}}_1$ and $x^{\{v,w\}}_{k-1}$, respectively. Let $c'$ be a total dominator edge coloring of  $G^{\frac{1}{k}}$. The two following cases can be occured: either the restriction of $c'$ to edges of $P^{\{v,w\}}$ is a total dominator edge coloring and so we have the result, or not. If the restriction of $c'$ to edges of $P^{\{v,w\}}$ is not a total dominator coloring then since $c'$ is a total dominator edge coloring  of $G^{\frac{1}{k}}$ we conclude that at least one of edges $vx^{\{v,w\}}_1$ and $wx^{\{v,w\}}_{k-1}$, as the edges of the induced subgraph $P^{\{v,w\}}$, are not adjacent to every edge of some color class. Without loss of generality we assume that the edge $vx^{\{v,w\}}_1$, as the edge of the induced subgraph $P^{\{v,w\}}$,  is not  adjacent to every vertex of some color class. But $c'$ is a total dominator coloring of $G^{\frac{1}{k}}$ so the edge $vx^{\{v,w\}}_1$ is adjacent to every edge of some color class, as the edge of  $G^{\frac{1}{k}}$. Hence there is a new color for an adjacent edge of $vx^{\{v,w\}}_1$, except the edge $x^{\{v,w\}}_1x^{\{v,w\}}_2$. Thus if we use this new color for the edge  $x^{\{v,w\}}_1x^{\{v,w\}}_2$ and consider the restriction of $c'$ for the remaining edges of superedge $P^{\{v,w\}}$, then $P^{\{v,w\}}$ has a total dominator edge coloring. Therefore the total edge coloring $c'$ has at least  $\chi'^t_d(P_{k+1}) $ colors.\qed
  
  \medskip
    	The lower bound of Theorem \ref{ubboun2p} is sharp for $P_2$ and by the following Proposition we  show that the upper bound of this Theorem is sharp for $G=K_{1,n}$ and $k=3$.

  \begin{proposition}
  	For every $n\geq 3$,  $\chi'^t_d(K_{1,n}^{\frac{1}{3}}) = 2n$.
  \end{proposition}
  
  \proof  Let $e_1, \ldots , e_n$ be the pendant edges of $K_{1,n}^{\frac{1}{3}}$. The adjacent edges to $e_i$ is denoted by $f_i$, and the adjacent edge to $f_i$  is denoted by $g_i$ for any $1\leq i \leq n$. Since  edge  $f_i$ is the only edge adjacent to $e_i$, so the color of $f_i$ should not be used for any other edges of graph, where $1\leq i \leq n$. Thus we color the edges $f_1, \ldots , f_n$ with colors $1, \ldots , n$, respectively, and do not use these colors any more. For every $1\leq i \leq n$, the edge $f_i$ is adjacent to $e_i$ and $g_i$, thus we need a new color for at least one of $e_i$ and $g_i$. So we need at least $2n$ color to have a TDE-coloring of $K_{1,n}^{\frac{1}{3}}$. Now for every $e_i$ and $g_i$ we use the new color $i+n$. Obviously this is a TDE-coloring of $K_{1,n}^{\frac{1}{3}}$ and we have the result. \qed
  
  Here we improve the  lower bound of Theorem \ref{ubboun2p} for  $k\geq 10$ .
  \begin{theorem}\label{lboundp}
  	If  $G$ is a connected graph with $m$ edges and maximum degree $\Delta(G)$ and $k\geq 10$, then 
  	\begin{equation*}
  	m(\chi'^t_d(P_{k-1})-2) + 2 \leq \chi'^t_d(G^{\frac{1}{k}}).
  	\end{equation*}
  \end{theorem}
  \proof Let $e=vw$ be an edge of $G$. We consider the superedge $P^{\{v,w\}}$ with vertex set $\{v, x^{\{v,w\}}_1, \ldots , x^{\{v,w\}}_{k-1}, w\}$. It is clear that $P^{\{v,w\}}\setminus \{v,w\}$ is the path graph $P_{k-1}$. Since we use  repetitious colors for the edges $x^{\{v,w\}}_1x^{\{v,w\}}_2$ and $x^{\{v,w\}}_{k-2}x^{\{v,w\}}_{k-1}$ in the  TDE-coloring  of paths, so we need at least $\chi'^t_d(P_{k-1})-2$ colors for each superedges and we cannot use these colors anymore. Also we need two colors for edges$x^{\{v,w\}}_1x^{\{v,w\}}_2$ and $x^{\{v,w\}}_{k-2}x^{\{v,w\}}_{k-1}$ and some other edges  hence  the result follows.\qed

  \begin{theorem}\label{lower1p}
  	If  $G$ is a connected graph with $m$ edges and maximum degree $\Delta(G)$ and $k\geq 10$, then  
  	\begin{equation*}
  	\chi'^t_d(G^{\frac{1}{k}})\geq \left\{
  	\begin{array}{ll}
  	m(\frac{k}{2})+2 & k\equiv 0 ~(\text{mod} ~4)\\
  	m(\frac{k-1}{2})+2 &k\equiv 1 ~(\text{mod} ~4)\\
  	m(\frac{k-2}{2})+2 & k\equiv 2 ~(\text{mod} ~4)\\
  	m(\frac{k-1}{2})+2 & k\equiv 3 ~(\text{mod} ~4).
  	\end{array}\right.
  	\end{equation*}
  \end{theorem}
  \proof 
  It follows by Theorems \ref{Pn} and \ref{lboundp} .\qed

  \begin{theorem}\label{uboundp}
  	If  $G$ is a connected graph with $m$ edges with maximum degree $\Delta(G)$ and $k\geq 10$  , then 
  	\begin{equation*}
  	\chi'^t_d(G^{\frac{1}{k}})\leq m(\chi'^t_d(P_{k+1})-2) +\Delta(G).
  	\end{equation*}
  \end{theorem}
  \proof As we see in the TDE-coloring  of paths, we can use the same color for the pendant edges. So we assign the colors $1,2,\ldots,\Delta(G)$ to all the edges incident  to the vertices belong to $G$ and we color other edges of any superedges with $\chi'^t_d(P_{k+1})-2$ colors. This is a TDE-coloring  for $G^{\frac{1}{k}}$ and hence  the result follows.\qed

  \begin{theorem}
  	If  $G$ is  a connected graph with $m$ edges and $k\geq 10$, then  
  	\begin{equation*}
  	\chi'^t_d(G^{\frac{1}{k}})\leq \left\{
  	\begin{array}{ll}
  	\frac{mk}{2}+\Delta(G), & k\equiv 0 ~(\text{mod} ~4)\\
  	m(\frac{k+1}{2})+\Delta(G), &k\equiv 1 ~(\text{mod} ~4)\\
  	m(\frac{k+2}{2})+\Delta(G), & k\equiv 2 ~(\text{mod} ~4)\\
  	m(\frac{k+1}{2})+\Delta(G), & k\equiv 3 ~(\text{mod} ~4).
  	\end{array}\right.
  	\end{equation*}
  \end{theorem}
  \proof 
  It follows by Theorems \ref{Pn} and \ref{uboundp} .\qed

  \begin{theorem}\label{frac}
  	For any $k\geq 4$, $\chi'^t_d(G^{\frac{1}{k}}) \leq \chi'^t_d(G^{\frac{1}{k+1}})$.
  	
  \end{theorem}
  
\begin{figure}[h]
  	\hspace{1.8cm}
  	\begin{minipage}{4.5cm}
  		\includegraphics[width=\textwidth]{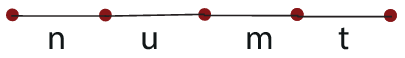}
  	\end{minipage}
  	\hspace{1cm}
  	\begin{minipage}{6cm}
  		\includegraphics[width=\textwidth]{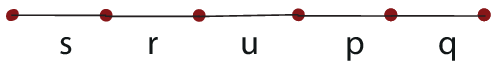}
  	\end{minipage}
  	\caption{\label{P5Conj} \small A part of a superedge in the proof of Theorem \ref{frac}.}
  \end{figure}

  \proof 
  First we give a TDE-coloring  to the edges of $G^{\frac{1}{k+1}}$. Let  $P^{\{v,w\}}$ be an arbitrary superedge of $G^{\frac{1}{k+1}}$ with vertex set  $\{v, x^{\{v,w\}}_1, \ldots , x^{\{v,w\}}_{k}, w\}$. We have the following cases:
  \begin{itemize}
  	\item[Case 1)]
  	There exists an edge $u\in \{ x^{\{v,w\}}_1x^{\{v,w\}}_2, \ldots , x^{\{v,w\}}_{k-1}x^{\{v,w\}}_{k}\}$ such that other edges of graph are not adjacent to all edges with color class of the edge $u$. Consider the graph in Figure \ref{P5Conj}. Suppose that the edge $u$ has the color $i$ and the edge $n$ has the color $\alpha$. The edge $m$ is adjacent to all edges with color class $j$ and $j\neq i$ and the edge  $n$ is adjacent to  all edges with color class $k$ and $k\neq i$. Since $k\geq 4$, without loss of generality, suppose that $m\neq vx^{\{v,w\}}_1$. We have two subcases:

  	\begin{itemize}
  		\item[Subcase i)]
  		The color of the edge $m$ is not $\alpha$. In this case, we  make $G/u$   and do not change the color of any edges. So without adding a new color we have a TDE-coloring for this new graph.
  		\item[Subcase ii)]
  		The color of the edge  $m$ is $\alpha$. Since the edge $u$ is adjacent to  color class $\alpha$, so any other edges does not  have color $\alpha$. In this case, by making $G/m$ and keeping the color of any edges as before, we have a TDE-coloring  for this new graph. Because the edge  $t$ is adjacent to color class which is not $\alpha$, the color of the edge  $t$ is not $i$ (because if  the color of the edge  $t$ is  $i$ it has contradiction with our assumptions), the edge $n$ is adjacent to  all edges with color class $k$ and  the edge $u$ is adjacent to all edges with color class $\alpha$.
  	\end{itemize}
  	
  	\item[Case 2)]
  	For every edge $u\in \{ x^{\{v,w\}}_1x^{\{v,w\}}_2, \ldots , x^{\{v,w\}}_{k-1}x^{\{v,w\}}_{k}\}$, there exists an edge such that is adjacent to all edges with color of edge $u$. Consider  the graph in Figure \ref{P5Conj}. Suppose that the edge $u$ has the color $i$ and the  edge $p$ has the color $j$ and the edge $p$ is adjacent to all edges with color $i$. We have two subcases:
  	\begin{itemize}
  		\item[Subcase i)]
  		The color of the edge $q$ is not $i$. We make $G/r$   and do not change the color of any edges. So without adding a new color we have a TDEC for this new graph since there is no other edges with color $i$.
  		\item[Subcase ii)]
  		The color of the edge $q$ is  $i$. In this case the edge $r$ is adjacent to color class of edge $s$ and the color of the edge $s$ does not use for other edges. Now we make $G/u$   and do not  change the color of any edges. Now we consider the color of edge $r$. If the color of $r$ is $j$, then we change it to $i$ and since obviously the edge $s$ was adjacent to a color class except $j$, so we have a TDE-coloring. If the color of the edge $r$ is not $j$ we do not  change the color of that and we have a TDE-coloring  again. 
  	\end{itemize}
  \end{itemize}
  Now we do the same algorithm for all superedges.  So we have $\chi'^t_d(G^{\frac{1}{k}}) \leq \chi'^t_d(G^{\frac{1}{k+1}})$.  \qed

\begin{figure}
 	\begin{center}
 		\includegraphics[width=1.5in]{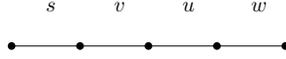}
 		\caption{A superedge in $G^{\frac{1}{4}}$.}
 		\label{P5e}
 	\end{center}
 \end{figure}

   \begin{theorem}
 For any graph $G$, 	$\chi'^t_d(G^{\frac{1}{3}}) \leq \chi'^t_d(G^{\frac{1}{4}}).$
  \end{theorem}
  \proof 
 First we give a TDE-coloring  to the edges of $G^{\frac{1}{4}}$. Let  $P^{\{i,j\}}$ be an arbitrary superedge of $G^{\frac{1}{4}}$ with edge set  $\{s,v,u, w\}$ (see Figure \ref{P5e}) and suppose that  the edge $v$ has the color $\alpha$. We have the following cases:
 \begin{itemize}
 	\item[Case 1)]
 The edges 	$u$ and $s$ are adjacent with an edge with a color class which is not $\alpha$. we have two subcases:
 	\begin{itemize}
 		\item[Subcase i)]
 		The color of edges $u$ and $s$ are different.  In this case, we make $G/v$ and don't change the color of any edges. So we have a TDE-coloring for this new graph. Because  two edges $u$ and $s$ are adjacent with an edge with  color class which is not $\alpha$.
 		\item[Subcase ii)]
 		The color of edges $u$ and $s$ are the same. Suppose that $u$ and $s$ have color $\beta$. In this case $\beta$ does not use for any other edges. So $w$ is adjacent with an edge with color class except $\beta$. Now we make $G/u$. So we have a TDE-coloring for this new graph.
 	\end{itemize}
 	
 	\item[Case 2)]
 The edge 	$u$ is adjacent to  all edges with color class $\alpha$. we have two subcases:
 	\begin{itemize}
 		\item[Subcase i)]
 		The color of the edge $w$ is not $\alpha$. Suppose that the edge $u$ has color $\gamma$.  If the edge $v$ is adjacent with all edges with color $\gamma$, and if the color of $s$ is  $\gamma$, we make $G/u$. But  if the color of edge $s$ is not $\gamma$, then we make $G/u$ and assign  the color $\gamma$ to the edge $w$. So we have a TDE-coloring for this new graph. If the edge $v$ is adjacent to  all edges with color except $\gamma$ (edge $s$), then we make $G/u$. So we have  a TDE-coloring for this new graph.
 		\item[Subcase ii)]
 		The color of the edge  $w$ is $\alpha$. We have two new cases. First, the edge $v$ is adjacent to an edge with color class $\gamma$. Any adjacent edge with $w$ is not adjacent to  edge with color class $\alpha$ (except $u$). So we make $G/u$ and assign  the color $\gamma$ to $w$. This is a TDE-coloring for this new graph.  Second, $v$ is not adjacent with color class $\gamma$. So the color of the edge $s$ does not use any more. Also the edge $s$ is not adjacent to edge  with color class $\alpha$. So we make $G/v$. This is a TDE-coloring for this new graph.
 	\end{itemize}
 	
 	\item[Case 3)]
 The edge	$s$ is adjacent to  all edges with color class $\alpha$. We have two subcases:
 	\begin{itemize}
 		\item[Subcase i)]
 		If $v$ is the only edge which has color $\alpha$,  then we make $G/u$ when $v$ is adjacent with color class of edge $s$ and make $G/s$ when $v$ is adjacent with color class of edge $u$. So this is a TDE-coloring for this new graph.
 		\item[Subcase ii)]
 		If there exist some edges with color $\alpha$, then the edge $u$ is adjacent with color class except $\alpha$. So we make $G/v$. This is a TDE-coloring for this new graph.
 	\end{itemize}
 \end{itemize}
 We apply this TDE-coloring for all superedges. So we obtain a TDE-coloring for  $G^{\frac{1}{3}}$. Therefore we have $\chi'^t_d(G^{\frac{1}{3}}) \leq \chi'^t_d(G^{\frac{1}{4}})$.  \qed

\begin{figure}
 	\begin{center}
 		\includegraphics[width=1.5in]{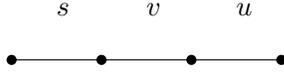}
 		\caption{A superedge in $G^{\frac{1}{3}}$.}
 		\label{P4e}
 	\end{center}
 \end{figure}

   \begin{theorem}
 For any graph $G$, 	$\chi'^t_d(G^{\frac{1}{2}}) \leq \chi'^t_d(G^{\frac{1}{3}}).$
  \end{theorem}
  \proof 
 First we give a TDE-coloring  to the edges of $G^{\frac{1}{3}}$. Let  $P^{\{i,j\}}$ be an arbitrary superedge of $G^{\frac{1}{3}}$ with edge set  $\{s,v,u\}$ (see Figure \ref{P4e}) and suppose that  the edge $v$ has the color $\alpha$. We have the following cases:
 \begin{itemize}
 	\item[Case 1)]
 The edges 	$u$ and $s$ are adjacent with an edge with a color class which is not $\alpha$. we have two subcases:
 	\begin{itemize}
 		\item[Subcase i)]
 		The color of edges $u$ and $s$ are different.  In this case, we make $G/v$ and don't change the color of any edges. So we have a TDE-coloring for this new graph. Because  two edges $u$ and $s$ are adjacent with an edge with  color class which is not $\alpha$.
 		\item[Subcase ii)]
 		The color of edges $u$ and $s$ are the same. Suppose that $u$ and $s$ have color $\beta$. In this case any other edges is not adjacent with color class $\beta$, because $i$ and $j$ are not adjacent vertices (Because of the definition of $G^{\frac{1}{3}}$). Now we make $G/u$. So we have a TDE-coloring for this new graph.
 	\end{itemize}
 	
 	\item[Case 2)]
The edge	$s$ is adjacent to  all edges with color class $\alpha$. We have two subcases:
 	\begin{itemize}
 		\item[Subcase i)]
 		If $v$ is the only edge which has color $\alpha$,  then we make $G/u$ when $v$ is adjacent with color class of edge $s$ and make $G/s$ when $v$ is adjacent with color class of edge $u$. So this is a TDE-coloring for this new graph.
 		\item[Subcase ii)]
 		If there exist some edges with color $\alpha$, then the edge $u$ is adjacent with color class except $\alpha$. If the edges $u$ and $s$ have the same color then we make $G/u$ and if $u$ and $s$ have different colors, then we make $G/v$. This is a TDE-coloring for this new graph.
 	\end{itemize}
 	
 	\item[Case 3)]
 The edges $u$ and	$s$ are adjacent to  all edges with color class $\alpha$. So there is no other edge with color $\alpha$. We have two subcases:
 	\begin{itemize}
 		\item[Subcase i)]
 		The edges $u$ and $s$ have the same color then we make $G/u$.
 		\item[Subcase ii)]
 		The edges $u$ and $s$ have different colors, then we make $G/u$ when $v$ is adjacent with color class of edge $s$ and make $G/s$ when $v$ is adjacent with color class of edge $u$.
 	\end{itemize}
 \end{itemize}
 We apply this TDE-coloring for all superedges. So we obtain a TDE-coloring for  $G^{\frac{1}{2}}$. Therefore we have $\chi'^t_d(G^{\frac{1}{2}}) \leq \chi'^t_d(G^{\frac{1}{3}})$.  \qed

\end{document}